\documentclass[11pt,reqno]{amsart}

\makeindex
\usepackage{amsmath}
\usepackage{amsthm}
\usepackage{amssymb}
\usepackage{amscd}
\usepackage{amsfonts}
\usepackage{amsbsy}
\usepackage[latin1]{inputenc}
\usepackage{enumerate}
\usepackage{indentfirst}
\usepackage{euscript}

\newtheorem {theorem} {Theorem}
\newtheorem {proposition} [theorem]{Proposition}

\newtheorem {definition} [theorem]{Definition}

\newtheorem {remark} [theorem]{Remark}

\newcommand{\R}{\ensuremath{\mathbb{R}}}
\newcommand{\T}{\ensuremath{\mathbb{T}}}

\def\a{{\alpha}}

\def\p{{\partial}}

\def\e{{\varepsilon}}

\def\prr#1{\frac{\partial }{\partial #1}}

\title[Lines of Curvature on Canal Surfaces ]
{Lines of Principal Curvature on Canal Surfaces in ${\mathbb
R^3}$}

 \author[R. Garcia, J. Llibre  and J. Sotomayor]
 {Ronaldo  Garcia, Jaume Llibre  and Jorge Sotomayor}


 \thanks{\noindent  E-mail: ragarcia@mat.ufg.br, sotp@ime.usp.br, jllibre@mat.uab.es\\
 .\hskip .3cm  AMS Classification(2000): 53C12, 34D30, 53A05, 37C75.}
\begin{document}
\maketitle

\begin{abstract}
In this paper are  determined the principal curvatures and
principal curvature lines on {\it canal surfaces} which are the
envelopes of families of spheres with variable radius and centers
moving along  a closed regular curve in $\mathbb R^3$. By means of
a connection of the differential equations for these curvature
lines and real Riccati equations, it is established that canal
surfaces have at most two isolated periodic principal lines.
Examples of canal surfaces with two simple and one double periodic
principal lines are given. \vskip .1cm

\end{abstract}

\section{Introduction}

The study of principal curvature lines, along which a surface in
$\mathbb R^3$ bends extremely and their umbilic singularities was
founded by Monge, Dupin and Darboux. See (Gray  1998) and
(Sotomayor  2003)
 for references.  For the basic facts about principal
curvature lines on surfaces the reader is addressed to (do Carmo
1976), (Spivak 1979) and  (Struik 1988).

\medskip

As a consequence of the work of Monge and Dupin the lines of
curvature on quadrics and toroidal --Dupin Cyclides-- surfaces
were determined. See (Fischer 1986, Chap. 3),  for an outline of
the theory and for a collection of remarkable illustrations.

\medskip

In (Gutierrez and Sotomayor 1982, 1991, 1998) ideas originating in
the Qualitative Theory of Differential Equations and Dynamical
Systems, such as Structural Stability and Genericity, were
incorporated into the subject; see also (Garcia and Sotomayor
2002). Historical comments, going up to some recent developments
on principal curvature lines, can be found in (Sotomayor 2003).

\medskip

The global dynamic complexity of principal curvature  in simple
smooth surfaces was illustrated in (Gutierrez and Sotomayor  1991,
1998).
 This includes examples of  recurrent (dense) principal
curvature lines on spheroidal surfaces, which are small
perturbations of ellipsoids of revolution and canal surfaces with
constant radial functions. The methods established in these works
show how to make hyperbolic a principal cycle by means of a small
smooth perturbation. This leads to smooth toroidal immersions with
arbitrary large number of principal isolated principal cycles.

\medskip

In this paper is improved and reproved a result that goes back to
Vessiot, establishing that principal curvature lines on canal
surfaces immersed in $\mathbb R^3$ (see Section \ref{sec:regcan},
Definition \ref{def:env}) verify a Riccati equation (see Section
\ref{sec:pc}, Remark \ref{rm:ric}). In fact, (Vessiot 1919) is the
first source known by the authors for the connection between
Riccati equations and principal curvature lines on canal surfaces.
The improvement consists in the formulation of precise conditions
for canal surfaces  be regular immersions and also
 to have umbilic points (see Section \ref{sec:pc}, Theorem
\ref{prop:envelope} and Remark \ref{rm:umb}). A consequence of the
Riccati structure for principal curvature lines on canal immersed
surfaces implies that the maximal number  of isolated periodic
principal lines is $2$. Examples of canal surfaces with two
(simple i.e. hyperbolic) and one (double i.e. semi--stable)
principal periodic lines, absent in (Vessiot 1919)  and also in
later references, are given in this paper (see Section
\ref{sec:ex}, Proposition \ref{prop:tt0}).

\section{ Regular Canal Surfaces }\label{sec:regcan}

Consider the space $\mathbb R^3$  endowed with the Euclidean inner
product $<\, , \,>$ and  norm $|\,\,\,| = <\, ,\,>^{1/2}$ as well
as with a canonical orientation.  The wedge product $\wedge$ of
vectors is defined relative to this orientation.

\medskip

Let ${\mathbf c}$ be a smooth regular closed curve immersed  in
$\mathbb R^3$, parametrized by arc length $s\in [0,L]$. This means
that
\begin{equation} \label{eq:ct}
{\mathbf c}^{\prime} (s)={\mathbf t}(s), \; |{\mathbf t}(s)|=1,
\;{\mathbf c}(L)={\mathbf c}(0)\;  .
\end{equation}
Assume also that the curve is {\it bi-regular}. That is:
\begin{equation} \label{eq:bireg}
\kappa(s) = |{\mathbf t}^{\prime} (s)| >0.
\end{equation}

Along ${\mathbf c}$ is defined  its  moving {\it Frenet frame}
$\{{\mathbf t},{\mathbf n},{\mathbf b}\}$. Following (Spivak 1979)
and (Struik 1988), this  frame is {\it positive}, {\it
orthonormal} and verifies {\it Frenet equations}:
\begin{equation} \label{eq:fr}
{\mathbf t}^{\prime} (s) = \kappa(s){\mathbf n}(s), \; {\mathbf
n}^{\prime} (s) = -\kappa(s){\mathbf t}(s) + \tau (s) {\mathbf
b}(s), \; {\mathbf b}^{\prime} (s) = -\tau (s){\mathbf n}(s).
\end{equation}

Equations  \eqref{eq:ct}  to   \eqref{eq:fr}  define the {\it unit
tangent}, ${\mathbf t}$, {\it principal normal}, ${\mathbf n}$,
{\it curvature},  $\kappa$,  {\it binormal}, ${\mathbf b}={\mathbf
t}\wedge {\mathbf n}$, and {\it torsion},  $\tau$,  of the
immersed curve ${\mathbf c}$.

\begin{proposition}\label{prop:immer}
Let $r(s)>0$ and $\theta(s) \in \; ]0,\pi[$ be smooth functions of
 period  $L$. The mapping $\alpha:\mathbb T^2 = \mathbb S^1\times
\mathbb S^1 \to\mathbb R^3 ,$ defined on $\mathbb R^2$ modulo
$L\times 2\pi$ by
\begin{equation}\label{eq:al}
\alpha(s, \varphi)={\mathbf c}(s)+ r(s)\cos\theta(s) {\mathbf
t}(s)+  r(s)\sin\theta(s)[\cos \varphi\, {\mathbf n}(s)+
\sin\varphi\, {\mathbf b}(s)],
\end{equation}
is tangent to the sphere of center ${\mathbf c}(s)$ and radius
$r(s)$ if and only if
\begin{equation}\label{eq:r1}
\cos\theta(s) = -r^{\prime} (s).
\end{equation}
Assuming \eqref{eq:r1}, with $r^{\prime} (s) < 1$, $\alpha$ is an
immersion provided
\begin{equation}\label{eq:reg1}
\kappa(s) <\frac{1-r'(s)^2 - r(s)r''(s)}{r(s)\sqrt{1-r'(s)^2}}.
\end{equation}
\end{proposition}

\begin{proof}
Calculations using equations  \eqref{eq:fr}  give
$$
\alpha_s= \frac{\p \alpha}{\p s} = \alpha_s^1{\mathbf t}(s)+
\alpha_s^2{\mathbf n}(s)+ \alpha_s^3{\mathbf b}(s),
$$
where
\begin{equation}\label{eq:alsteta}
\aligned \alpha_s^1=\, &1+r^{\prime} \cos\theta-r \theta^{\prime}
\sin\theta-r \kappa \sin\theta \cos\varphi,\\
\alpha_s^2=\, &r^{\prime}\sin\theta \cos\varphi+r \theta^{\prime}
\cos\theta\cos\varphi -r\tau \sin\theta \sin\varphi +r\kappa
\cos\theta ,\\
\alpha_s^3=\, &r^{\prime}\sin\theta \sin\varphi+r\theta^{\prime}
\cos\theta \sin\varphi+r \tau\sin\theta \cos\varphi.
\endaligned
\end{equation}
Here, of course, $r$, $\theta$, $\tau$, $\kappa$ are functions of
$s$. Also,
\[
\alpha_\varphi = \frac{\p \a}{\p \varphi} =
r(s)\sin\theta(s)[-\sin\varphi \,{ \mathbf n}(s)+ \cos\varphi{
\mathbf b}(s)].
\]

The unit  normal vector  pointing  inward the sphere $|p-{ \mathbf
c}(s)|=r(s)\,$ with center ${ \mathbf c}(s)$ and radius $r(s)$ at
$p= \alpha(s,\varphi)$ is $N_\alpha =({ \mathbf
c}(s)-\alpha(s,\varphi))/r(s)$, which is given by:
\begin{equation}\label{eq:normal}
N_\alpha = -\cos\theta(s)\, {\mathbf t}(s) - \sin\theta(s)
[\cos\varphi  \,{\mathbf n}(s) + \sin\varphi  \,{ \mathbf b}(s)].
\end{equation}
Clearly, $<\alpha_{\varphi},N_\alpha> =0$.  Calculation gives
$<\alpha_{s},N_\alpha> =-(\cos\theta + r^{\prime}) $. Therefore
the condition of tangency at $(s,\varphi)$  of $\alpha$ to the
sphere $|p-{\mathbf c}(s)|=r(s)\,$ is $\cos\theta(s) =
-r^{\prime}(s)$.

\medskip

Additional calculation gives:
\begin{equation}\label{eq:FG}
\aligned
 F(s,\varphi) &= <\alpha_s,\alpha_\varphi>= \tau
r^2\sin^2\theta+\kappa(s)r^\prime r^2\sin\theta\sin\varphi, \\
G(s,\varphi) &=<\alpha_\varphi,\alpha_\varphi>=
r^2\sin^2\theta=r^2(1-{r^\prime}^2).
\endaligned
\end{equation}

At this point it is appropriate to write the mapping $\alpha$ and
its derivatives involving only functions  of $ r(s)$ rather than
of $\theta(s)$,  replacing the expressions $\cos\theta(s)
=-r'(s),\, \sin\theta(s) = (1-(r')^2)^{1/2}$ in equations
\eqref{eq:al}  and  \eqref{eq:alsteta}.

\medskip

The expression for $E$ is as follows:
\begin{equation}\label{eq:E}
\aligned
E(s,\varphi) =  & <\alpha_s,\alpha_s> = \kappa^2 r^2 (1-{r^\prime}^2) \cos^2\varphi \, +\\
& 2\, \kappa  r [r r^\prime  -(1-{r^\prime}^2)/\sqrt{1-{r^\prime }^2}]\cos\varphi \, +\\
& 2\kappa  \tau  r^2  r^\prime \sqrt{1-{r^\prime }^2}\sin\varphi\
+(1-{r^\prime}^2) + 2rr^{\prime\prime} +\\
& \frac{ r^2}{(1-{r^\prime}^2)}[\tau^2 (1-{r^\prime}^2)^2 +
{r^\prime}^2 \kappa^2 (1-{r^\prime}^2) +{r^{\prime\prime}}^2].
\endaligned
\end{equation}

Additional calculation and simplification using equations
 \eqref{eq:FG}, with the tangency condition  \eqref{eq:r1}
imposed, and  \eqref{eq:E}  gives
\begin{equation}\label{eq:area}
EG-F^2= r(s)^4[\sqrt{1-r'(s)^2} \kappa(s) \cos\varphi + r''(s)-(1-
r'(s)^2)/r(s)]^2.
\end{equation}
This expression, which is equal to $|\partial{\alpha}/\partial s
\wedge \partial{\alpha}/\partial \varphi|^2$, vanishes if and only
if
\begin{equation}\label{eq:rt1}
\cos\varphi = \frac{1-r'(s)^2 -r(s)r''(s)}{\sqrt{1-r'(s)^2}
\kappa(s)r(s)}.
\end{equation}

Condition  \eqref{eq:reg1}  states that the absolute value of the
right--hand member of  \eqref{eq:rt1}  is larger than $1$, which
by
 \eqref{eq:area}  implies the linear independence of $\alpha_s$
and $\alpha_\varphi$ at every point $(s,\varphi)$, namely that
$\alpha$ is an immersion.
\end{proof}

\begin{definition} \label{def:env}
A mapping such as $\alpha$, of $\,\, \T^2$ into $\R^3$, satisfying
conditions  \eqref{eq:r1} and \eqref{eq:reg1}  will be called an
{\it immersed canal surface} with {\it center along} ${\mathbf
c}(s)$ and {\it radial function} $r(s)$. When $r$ is constant, it
is called an {\it immersed tube}. Due to the tangency condition
\eqref{eq:r1}, the  {\it immersed canal surface} $\alpha$ is the
{\it envelope} of the family of spheres of radius $r(s)$ whose
centers range along the curve  ${\mathbf c}(s)$.
\end{definition}

\begin{remark}\label{rm:ve}
Proposition \ref{prop:immer} is partially found in (Vessiot 1919).
There, however, the regularity condition \eqref{eq:reg1} was
overlooked. Other references for canal surfaces are (Blaschke
1929)  and (Gray 1998).
\end{remark}

\section{Principal Lines  on Regular Canal Surfaces}
\label{sec:pc}

In this work the positive orientation on the torus $\mathbb T^2 =
\mathbb S^1\times \mathbb S^1$ is defined by the ordered tangent
frame $\{\partial/\partial s ,\partial/\partial \varphi\}$.
Therefore the {\it positive unit normal} --or Gaussian map--
$N_{\alpha}$  of the immersion $\alpha$ is defined by
$|\partial{\alpha}/\partial s \wedge \partial{\alpha}/\partial
\varphi|N_{\alpha}= \partial{\alpha}/\partial s \wedge
\partial{\alpha}/\partial \varphi$. By the tangency condition
 \eqref{eq:r1}  this unit vector is given by \eqref{eq:normal},
which can be written as follows:
\begin{equation}\label{eq:normalr}
N_\alpha = r'(s)\, {\mathbf t}(s) - (1-(r'(s))^2)^{1/2}
[\cos\varphi\, {\mathbf n}(s)+ \sin\varphi\,  {\mathbf b}(s)].
\end{equation}
It points inwards the toroidal surface defined by $\alpha$.

\medskip

Below will be studied the global behavior of the principal
curvature lines of $\alpha$.

\begin{theorem}\label{prop:envelope}
Let $\alpha:\mathbb  S^1\times \mathbb  S^1 \to\mathbb R^3$ be a
smooth immersion expressed  by \eqref{eq:al}. Assume the
regularity conditions \eqref{eq:r1} and \eqref{eq:reg1} as in
Proposition \ref{prop:immer} and also that

\begin{equation} \label{eq:nu}
k(s)<\Big| \frac{1-r'(s)^2-2r(s)r^{\prime\prime}(s) }{2 r(s)
\sqrt{ 1-r'(s)^2}}\Big|,
\end{equation}
The maximal principal curvature lines are the circles tangent to
$\partial/\partial \varphi$. The maximal principal curvature is
\[
k_2 (s) = 1/r(s).
\]
The minimal principal curvature lines are the curves tangent to
\begin{equation}\label{eq:clc}
V(s,\varphi)= \frac{\partial}{\partial s} - \left(\tau(s)+
\frac{r'(s)}{(1-r'(s)^2)^{1/2}}\kappa(s)\sin \varphi\right)
\frac{\partial}{\partial \varphi}.
\end{equation}
The expression
\begin{equation}\label{eq:d}
r(s)(1-r'(s)^2)^{1/2}\kappa(s)\cos\varphi
+r(s)r^{\prime\prime}(s)- (1-r'(s)^2)
\end{equation}
is negative, and the minimal  principal curvature is given by

\begin{equation}\label{eq:k1}
k_1 (s,\varphi) = \frac{\kappa(s)(1-r'(s)^2)^{1/2}\cos\varphi
+r^{\prime\prime}(s)}{r(s)(1-r'(s)^2)^{1/2}\kappa(s)\cos\varphi
+r(s)r^{\prime\prime}(s) - (1-r'(s)^2)}.
\end{equation}

There are  no umbilic points for $\alpha$: $k_1 (s,\varphi) < k_2
(s). $

\end{theorem}

\begin{proof}
Direct calculation gives:
$$
\frac{\partial N_\alpha}{\partial\varphi } + k_2(s)\frac {\partial
\alpha}{\partial \varphi} =0.
$$
By {\it Rodrigues equation}, see (do Carmo  1976) and (Struik
1988), the circles $s=$ constant are principal curvature lines of
$\alpha$, with principal curvature $k_2(s)$.

\medskip

Denoting the metric of $\alpha$  by $Eds^2+2Fdsd\varphi+
Gd\varphi^2$, it follows that the direction orthogonal to
$\partial/\partial \varphi$, giving the other principal direction,
is defined by the vector field
$$
G\prr s-F\prr \varphi,
$$
which by equation  \eqref{eq:FG}   is collinear with $
V(s,\varphi)$ in  \eqref{eq:clc}.

\medskip

Differentiation of  \eqref{eq:normalr}   using Frenet equations
 \eqref{eq:fr}  gives
\begin{equation}\label{eq:ns}
\aligned  \frac{\partial N_\alpha}{\partial s} = \, &
[r''+\kappa(1-{r'}^2)^{1/2})\cos\varphi]{\mathbf t}(s) +\\
&[(r'r''\cos \varphi  +\tau  (1-{r'}^2)\sin \varphi +\\
& r' \kappa (1-{r'}^2)^{1/2})/(1-{r'}^2)^{1/2}]{\mathbf n}(s) +\\
&[(-\tau  (1-{r'}^2)^{1/2}\cos\varphi+ r'
r''\sin\varphi)/(1-{r'}^2)^{1/2}]{\mathbf b}(s),
\endaligned
\end{equation}
and
$$
\frac{\partial N_\alpha}{\partial \varphi} = \sqrt{1-r'(s)^2}\,
\left[\sin\varphi\,{\mathbf n}(s)-\cos\varphi\,{\mathbf
b}(s)\right].
$$

Substitution, taking into account equation \eqref{eq:clc}, leads
to

\begin{equation}\label{eq:nv}
\frac{\partial N_\alpha}{\partial V}(s,\varphi)=\frac{\partial
N_\alpha}{\partial s} - \left[\tau +\frac{r'
}{(1-{r'}^2)^{1/2}}\kappa \sin \varphi\right]\frac{\partial
N_\alpha}{\partial \varphi}.
\end{equation}

 Similarly for $\alpha$,  calculation of its derivative $\frac{\partial  \alpha}{\partial
 V}$,
 gives

\begin{equation}\label{eq:av} \aligned
\frac{\partial  \alpha}{\partial V}(s,\varphi)=& A(s,\varphi)[
{\mathbf t}(s) + \frac{r^\prime(s)}{ (1-{r^\prime(s)}^2)^{1/2}}
(\cos\varphi
{\mathbf n}(s)+ \sin\varphi {\mathbf b}(s))],\\
A(s,v)=&  1-{r^\prime(s)}^2-r(s) r^{\prime\prime}(s) -\kappa(s)
r(s)\cos\varphi (1-{r^\prime(s)}^2)^{1/2}  \endaligned
\end{equation}

 Further calculation using   equations \eqref{eq:ns}, \eqref{eq:nv},
\eqref{eq:av} and the expression
 for $k_1(s,\varphi)$ given in
equation \eqref{eq:k1}, follows that
$$\frac{\partial N_\alpha}{\partial
V}(s,\varphi) + k_1(s,\varphi) \frac{\partial \alpha}{\partial
V}(s,\varphi) = 0.
$$
\noindent This is Rodrigues equation,  which establishes that, in
fact,  $k_1$ is a principal curvature.

  From  \eqref{eq:reg1}  it follows that  \eqref{eq:d}, the
denominator of $k_1$, is always negative.
 This implies that $k_1 < 1/r(s) =
k_2$.  Otherwise, $k_1 \geq k_2$ and, after direct manipulation,
by \eqref{eq:nu} this would lead to
\begin{equation} \label{eq:u}
\Big|\cos\varphi\Big| \geq \Big|
\frac{1-r'(s)^2-2r(s)r^{\prime\prime}(s) }{2\kappa(s) r(s) \sqrt{
1-r'(s)^2}}\Big|>1.
\end{equation}
This prevents the existence of
  umbilic points and justifies the names {\it
maximal}, for subscript $2$, and {\it minimal}, for subscript $1$,
given in the statement.
\end{proof}

\begin{remark}[Umbilic Points in Canal Immersions]\label{rm:umb}
The calculation of $k_1$ and the condition for the appearance of
umbilic points for regular canal surfaces has not been considered
in previous works on the subject. The discussion leading to
\eqref{eq:u} leads to
 the following equation for umbilic points:
\[
\cos\varphi =  \frac{1-r'(s)^2-2r(s)r^{\prime\prime}(s)
}{2\kappa(s) r(s)\sqrt{1-r'(s)^2}}.
\]
It gives a non--empty curve if condition \eqref{eq:nu} is not
imposed. This curve however consists of removable singularities
for the principal line fields which have smooth extensions to the
whole torus, given by $\partial/\partial \varphi$ and $V$. Under
generic conditions on the radial function, $r(s)$, and curvature,
$\kappa(s),$ functions it will be expected to appear curves of
umbilic points as those studied in Proposition 2 of  (Garcia and
Sotomayor 2005).
\end{remark}

\begin{remark}[Minimal Principal Foliation in terms of Differential
Forms] Written as a differential form, the vector field
\eqref{eq:clc} becomes:
\begin{equation}\label{eq:clf}
\omega_1= d\varphi+ \left[\tau(s) +\frac{r'(s)}
{(1-r'(s)^2)^{1/2}} \kappa(s)\sin \varphi\right]ds.
\end{equation}
When it is used another parameter $t$, say $T$--{\it periodic},
related with the arc length $s $ of ${\bf c}$  in the form
$s=s(t)$, with $\dot s = <\dot {\mathbf c}, \dot {\mathbf
c}>^{1/2} $,  the form in \eqref{eq:clf} writes as:
\begin{equation}\label{eq:clt}
\omega_1= d\varphi+ \left[\tau(t) + \frac{\dot r(t)}{(\dot s
(t)^2-\dot r (t)^2)^{1/2}}\kappa(t)\sin\varphi\right]\dot s(t) dt.
\end{equation}
\end{remark}

\begin{remark}[Vessiot]\label{rm:ric}
By the change of coordinates $\tan(\varphi/2) =z$, equation
\eqref{eq:clf} is transformed into the $L$--periodic Riccati
equation
\begin{equation}\label{eq:rtp}
z^\prime = -\frac 12 \tau(s)(1+z^2)+ \cot\theta(s)\kappa(s)z.
\end{equation}
Therefore, the solutions of equation \eqref{eq:clc} can be
obtained from the solutions of equation \eqref{eq:rtp} contained
in $[0,T]\times [-\pi, \pi]$.
\end{remark}

\section{Canal Surfaces with one and two Principal Cycles}\label{sec:ex}

In this section is carried out a discussion on the qualitative
properties of equation  \eqref{eq:clc}, or of its equivalent form
\eqref{eq:clt}, absent in (Vessiot 1919).

\medskip

Being equivalent to a periodic Riccati equation, this equation has
a {\it M\"{o}ebius  transformation} as {\it return map}  (see
Hille 1976) and therefore can have either:
\begin{itemize}
\item[(a)] all its solutions periodic, \item[(b)] all its
solutions dense,
 \item[(c)] two hyperbolic (simple)  periodic
solutions or
 \item[(d)] one semi--hyperbolic (double or semi--stable) periodic
solution.
\end{itemize}

An example of situation (a) is exhibited by the standard torus of
revolution. An example of (b) is  given in (Gutierrez and
Sotomayor 1991, 1998), for a canal surface of constant radius (a
{\it tube}) around  a  curve that is not bi-regular. Below will be
given examples of cases (c) and (d). This will also provide
examples of case (b) for bi-regular curves. The cost of this is
heavier calculation which, nevertheless,  is easy to corroborate
with Computer Algebra.
\medskip

The example consists in a deformation of the $T=2\pi$--{\it
periodic} plane elliptic curve
$$
{\mathbf c}(t)=(2\cos t, \sin t, 0),
$$
whose curvature $\kappa (t)=\kappa (t,0)$ is $\kappa (t) = 2/(4
-3\cos ^2t)^{3/2}.$

\begin{proposition}\label{prop:tt0}
Consider the three parameter family of canal surfaces
$S_{\e,\rho,\mu}$ around the curve  ${\mathbf c}_{\e}(t) = (2\cos
t, \sin t, \e \dot \kappa  (t))$, with radial function
$$
r(t,\mu)=\rho + \mu\, \dot \kappa (t), \quad \dot \kappa(t) =
-18\cos t\, \sin t /(4 -3\cos^2t)^{5/2}.
$$
There are two smooth curves $\e= \e_1 (\mu)= -\mu + O_1(\mu^2)$
and $\e= \e_2 (\mu)= \mu + O_2(\mu^2)$, such that for any $\rho,
\, \mu$ small and positive,  the canal surface $S_{\e,\rho , \mu}$
has two hyperbolic principal cycles for $\e\in ]\e_1(\mu),
\e_2(\mu)[$ and has one double principal cycle along the curves
$\e= \e_1 (\mu)$ and $\e= \e_2 (\mu)$.
\end{proposition}

\begin{proof} Standard calculation gives:
$$
\dot s (t,\e)= (4\sin^2t + \cos^2 t + (\e \ddot{ \kappa
}(t))^2)^{1/2} = (-3\cos^2 t+4)^{1/2} +O(\e ^{2}),
$$
for the element of arc length of ${\bf c}_{\e}$.

Also, the curvature of ${\mathbf c}_{\e}$ is
\begin{eqnarray*}
\kappa(t,\e) &=& \kappa(t,0)+ 648\e^2  \left[ \frac{
1458\cos^{12}t+486\cos^{10}t-4671\cos^8t}{(4 -
3\,\mathrm{cos}^{2}t)^{19/2}} \right. +\\
& &\left. \frac{432\cos^6t+4296\cos^4t-2016\cos^2t+16 }{(4 -
3\,\mathrm{cos}^{2}t)^{19/2}} \right]+ O(\e^3).
\end{eqnarray*}
Here will be needed only $\dot \kappa(t,0) =\dot \kappa(t)$ as
given above and  $\ddot \kappa(t)$ in the integral
\begin{equation}\label{eq:kint}
\int_0^{2\pi}\dot s (t,0)\ddot{\kappa}(t)\kappa(t)dt =
\int_0^{2\pi} \frac{36(9\cos^4t-4\cos^2t-4)}{(-3\cos^2t+4)^5}\, dt
= - \frac{8829}{2048}\pi.
\end{equation}

The torsion of ${\mathbf c}_{\e}$ is
$$
\tau (t,\e) = \tau(t,0)+\e\frac{ \partial \tau}{\partial \e
}(t,0)+ O(\e^2).
$$
Here will only be needed the expressions for $\tau(t,0)=0$ and $
\tau_1 (t) = \partial \tau/\partial \e (t,0).$

From direct calculations follows that
\begin{eqnarray}
\tau_1(t) &=& \frac{108\,(54\,\mathrm{cos}^{8}t +
207\,\mathrm{cos}^{6}t - 369\,\mathrm{cos}^{4} t+
68\,\mathrm{cos}^{2}t + 44)} {  (4 -
3\,\mathrm{cos}^{2}t)^{11/2}}, \nonumber\\
C &=& -\int_0^{2\pi} \dot s (t,0)\tau_1(t)dt  =
\frac{8829}{2048}\pi. \label{eq:ti}
\end{eqnarray}

The differential equation corresponding to  \eqref{eq:clt}  for
the $2\pi$--periodic canal surface $S_{\e,\rho,\mu}$ is
\begin{equation}\label{eq:ve}
\frac{d\varphi}{dt} =  W(t,\varphi, \e, \mu)\,:=  \dot s
(t,\e)\left(- \tau (t,\e)  - \frac{\mu\ddot \kappa(t) \kappa(t,\e)
\sin\varphi}{\sqrt{\dot s (t)^2 - \mu ^2 \ddot \kappa(t)
^2}}\right).
\end{equation}

Denote by $\Pi (\varphi_0, \e, \mu)$ the return map for equation
\eqref{eq:ve}. It is given by  $\Pi (\varphi_0, \e, \mu)=
\Phi(2\pi,\varphi_0, \e, \mu)$, where $\Phi(t,\varphi_0, \e, \mu)$
is the solution of \eqref{eq:ve} such that $\Phi(0,\varphi_0, \e,
\mu)= \varphi_0$. Therefore, the $2\pi$--periodic solutions are
given by the implicit surface $\Pi (\varphi_0, \e, \mu)-\varphi_0
=0$. Since $\Pi (\varphi_0, 0, 0) =\varphi_0$,  the Fundamental
Theorem of Calculus (Hadamard's Formula) implies that
\begin{equation}\label{eq:fold1}
\Pi (\varphi_0, \e, \mu)-\varphi_0= \e\int_0^1 \frac{\partial
\Pi}{\partial \e } (\varphi_0, u \e, u\mu)du +
\mu\int_0^1\frac{\partial \Pi }{\partial \mu}(\varphi_0, u \e,
u\mu) du.
\end{equation}

Write
\begin{eqnarray}
\frac{\partial \Pi }{\partial \e }(\varphi_0, u \e, u\mu) &=&
\frac{\partial \Pi}{\partial \e } (\varphi_0, 0, 0) +
P(\varphi_0,u,\e,\mu)\e + Q(\varphi_0,u,\e,\mu)\mu,
\label{eq:folde}\\
\frac{\partial \Pi } {\partial \mu }(\varphi_0, u \e, u\mu) &=&
\frac{\partial \Pi}{\partial \mu } (\varphi_0, 0, 0) +
R(\varphi_0,u,\e,\mu)\e + S(\varphi_0,u,\e,\mu)\mu.
\label{eq:foldm}
\end{eqnarray}

The expressions for the derivatives of the solutions of the
differential equations with respect to parameters applied to
 \eqref{eq:ve}, which in the present case, following classical
differential equations results, are the integrals of the
non--homogeneous linear
---called {\it variational}---  equations.  According
to (Sotomayor 1979, page 42), (Chicone 1999, page 337) or
(Coddington and Levinson 1955, page 30), where these classical
results are proved, these equations and their initial conditions
are as follows:
\begin{eqnarray*}
\dot\Phi_{\e}(t,\varphi_0,0,0) &=& W_{\varphi}(t,\varphi_0,0,0)
\Phi_{\e}(t,\varphi_0,0,0)+ W_{\e}(t,\varphi_0,0,0)\\
&=& -\dot s(t,0)\tau_1(t),\;\; \Phi_{\e}(0,\varphi_0,0,0) = 0, \\
\dot\Phi_{\mu}(t,\varphi_0,0,0) &=& W_{\varphi}(t,\varphi_0,0,0)
\Phi_{\mu}(t,\varphi_0,0,0) + W_{\mu}(t,\varphi_0,0,0)\\
&=& -\dot s (t,0)\ddot{\kappa}(t)\kappa(t)sin\varphi_0 , \;\;
\Phi_{\mu}(0,\varphi_0,0,0)= 0.
\end{eqnarray*}

Integrating these equations, taking into account equations
\eqref{eq:kint} and   \eqref{eq:ti}, leads to
\begin{eqnarray}
\frac{\partial \Pi}{\partial \e } (\varphi_0, 0, 0) &=&
\Phi_{\e}(2\pi,\varphi_0,0,0)= - \int_0^{2\pi} \tau_1(v)dv = C,
\label{eq:folde1}\\
\frac{\partial \Pi}{\partial \mu } (\varphi_0, 0, 0) &=&
\Phi_{\mu}(2\pi,\varphi_0,0,0) = C\sin\varphi_0. \label{eq:foldm1}
\end{eqnarray}
Writing $\e = \nu \mu$ and substituting into \eqref{eq:fold1},
taking into consideration the expressions \eqref{eq:folde1} and
\eqref{eq:foldm1}, obtain
$$
\Pi (\varphi_0, \e, \mu)-\varphi_0 = C \mu( \nu  +  \sin\varphi_0
+\nu Z(\nu,\mu, \varphi_0)).
$$
Therefore,  the surface of $2\pi$--periodic orbits, in the $(\nu,
\mu,\varphi)$--space consists on a plane $\mu=0$ crossing
transversally a regular sheet $\mathbb P$  at the sinusoidal curve
of equation $\nu+ \sin\varphi_0 =0$. The critical values of the
projection on the $(\nu, \mu) -$plane  of the surface $\mathbb P$,
coming from double (semi--stable)  periodic solutions, the {\it
fold} curve, must cross transversally the $\nu$--axis at points
corresponding to $\nu = \pm 1$,
 which are
the critical values of the projection of the sinusoidal curve.
Going back to the coordinates $(\e,\mu)$ gives the conclusion
formulated in the proposition.
\end{proof}

\begin{remark}[Arnold Tongues]
Due to equations \eqref{eq:folde1} and \eqref{eq:ve}, for $\rho$
small, the return map on the tube $S_{\e,\rho ,0}$ is a rotation
with its {\it rotation number} changing monotonically with $\e$
(taken small). Therefore it takes irrational values and the tube
presents case (b). It also takes rational values, $q/p$,  which
are the vertices of the {\it Arnold Tongues} (in the
$(\e,\mu)$--plane), corresponding to canal surfaces with periodic
closed principal lines winding $p$--times around the parallels
($t$--circles) and $q$--times around the meridians
($\varphi$--circles) of the torus. The sector with vertex at
$(0,0)$ established in Proposition \ref{prop:tt0} is the principal
tongue, with rotation number $0$. See (Chicone 1999, page 372).
\end{remark}

\vskip .5cm
 \centerline{\sc Acknowledgments}
\vskip .3cm
 The first and third authors are fellows of CNPq and are
partially supported by CNPq Grants PADCT 620029/2004-8 and
47.3824/2004-3. The second author is partially supported by a
DGICYT grant number MTM2005-06098-C02-01 and by a CICYT grant
number 2001SGR00173. The first  and second authors are also
partially supported by the joint project CAPES-MECD grants
071/04 and HBP2003-0017, respectively. \vskip .5cm

\section{References}

 \noindent {\sc  Blaschke W}. 1929.
{ Differential Geometrie, III, Differentialgeomtrie der Kreise und
Kugeln}, Springer--Verlag, Berlin.

\vskip .2cm

 \noindent  {\sc   Chicone C.}  1999. { Ordinary  Differential Equations with
Applications}, Springer-Verlag, New York. \vskip .2cm

 \noindent  {\sc Coddington E and  Levinson N.}  1955.
{ Theory of Ordinary Differential Equations}, Mc-Graw - Hill, New
York.

\vskip .2cm

 \noindent  {\sc  do Carmo M.} 1976.
{  Differential Geometry of Curves and Surfaces}, Prentice-Hall,
Englewood-Cliffs.

\vskip .2cm

 \noindent  {\sc   Fischer G.}
1986. { Mathematical Models},  Friedr. Vieweg $\&$ Sohn,
Braunschweig.

\vskip .2cm

 \noindent  {\sc   Gray A.}  1998.
{ Modern Differential Geometry of Curves and Surfaces with
Mathematica}, Second Edition, CRC Press, New York, 1998.

\vskip .2cm

 \noindent  {\sc   Garcia R and  Sotomayor J.}  2002.
{Lecture Notes on Differential Equations of Classical Geometry},
Preprint of Lecture Course delivered in São Carlos, Brazil.

\vskip .2cm
 \noindent  {\sc Garcia R and   Sotomayor J.}   2005.
 {\it On the patterns
of principal curvature around a curve of umbilic points}, An.
Acad. Brasil. Ci\^encias  { {\bf 77}:} 13--24.

\vskip .2cm
 \noindent   {\sc  Gutierrez C and  Sotomayor J.}  1982.
{\it Structural Stable Configurations of Lines of Principal
Curvature}, Asterisque { {\bf 98--99}:}  185--215.

\vskip .2cm
  \noindent  {\sc  Gutierrez C and   Sotomayor J.}  1991.
  {Lines of
Curvature and Umbilic Points on Surfaces}, in the $18$--th
Brazilian Math. Colloquium, Rio de Janeiro, IMPA. Reprinted in
1998, with update, as { Structurally Configurations of Lines of
Curvature and Umbilic Points on Surfaces}, Lima, Monografias del
IMCA.

\vskip .2cm
 \noindent  {\sc Hille E.} 1976.  { Ordinary Differential Equations in
the Complex Domain}, Dover Publications, New York.

\vskip .2cm
 \noindent  {\sc Sotomayor J.} 1979.  { Li\c c\~oes de Equa\c c\~oes Diferenciais Ordin\'arias},
 IMPA-CNPq, Projeto Euclides, Rio de Janeiro.

\vskip .2cm
  \noindent  {\sc   Sotomayor J.}   2003. {\it Historical Comments on
Monge's Ellipsoid and the Configuration of Lines of Curvature on
Surfaces Immersed in ${\mathbb R}^3$}. Mathematics ArXiv.\\
http://front.math.ucdavis.edu/math.HO/0411403.

\vskip .2cm
  \noindent  {\sc   Spivak M.}  1979. { Introduction to Comprehensive Differential Geometry},
Vol. II, III, Publish or Perish, Berkeley.

\vskip .2cm
 \noindent   {\sc  Struik D.}  1988. { Lectures on Classical
Differential Geometry}, Addison Wesley Pub. Co., Reprinted by
Dover Publications, Inc.

\vskip .2cm
  \noindent  {\sc   Vessiot E.}   1919.
   { Le\c cons de G\'eom\'etrie
Sup\'erieure}, Librarie Scientifique J. Hermann, Paris.

\newpage

\author{\noindent Ronaldo Garcia\\Instituto de Matem\'atica e Estat\'{\i}stica \\
Universidade Federal de Goi\'as, \\
CEP 74001--970, Caixa Postal 131 \\
Goi\^ania, Goi\'as, Brazil}

\vskip 0.5cm

\author{\noindent Jaume Llibre\\ Departament de Matem\`{a}tiques \\ Universitat
Aut\`{o}noma de Barcelona\\ 08193 Bellaterra, Barcelona, Spain}

\vskip 0.5cm
\author{\noindent Jorge Sotomayor\\ Instituto de Matem\'atica e Estat\'{\i}stica \\
Universidade  de S\~ao Paulo,\\
 Rua do Mat\~ao  1010,\\
Cidade Univerit\'aria, CEP 05508-090,\\
S\~ao Paulo, S. P, Brazil}

\end{document}